\documentclass[reqno]{amsart}

\usepackage[english]{babel}
\usepackage[numbers]{natbib}
\usepackage{bigints}
\usepackage{xcolor}
\usepackage{enumitem}
\usepackage{empheq}
\usepackage{array}
\usepackage{graphicx}
\usepackage{lipsum}
\usepackage{amssymb}
\usepackage{amsthm}
\usepackage{needspace}
\usepackage[makeroom]{cancel}
\newenvironment{myproof}[2] {\paragraph{\textbf{Proof of {#1} {#2} :}}}{\hfill$\square$}
\newenvironment{myproof1}[2] {\paragraph{\textbf{Proof of completeness of Example{#1} {#2} :}}}{\hfill$\square$}
\newtheorem{theorem}{Theorem}[section]

\newtheorem{question-non}[]{}

\newtheorem{observation}[theorem]{Remark}
\newtheorem{example}[theorem]{Example}

\usepackage{geometry}
\title{Invariant solutions of gradient $k$-Yamabe solitons}

\author{Tokura, W. $^{{1},\ast}$}
\address{$^{1}$ Universidade Estadual de Mato Grosso do Sul, 79150-000, Av. João Pedro Fernandes, 2101 - Centro, Maracaju, MS, Brazil.}
\email{williamisaotokura@hotmail.com $^{1}$}

\author{Barboza, M. $^{2}$}
\address{$^{2}$ Insituto Federal Goiano, 75790-000, Rodovia Geraldo Silva Nascimento Km 2,5, Uruta\'i, GO, Brazil.}
\email{marcelo.barboza@ifgoiano.edu.br $^{2}$}

\author{Batista, E. $^{3}$}
\address{$^{3}$ Universidade Federal de Goi\'as, IME, 131, 74001-970, Goi\^ania, GO, Brazil.}
\email{edbatista@gmail.com.br $^{3}$}

\author{Kai, P. $^{4}$}
\address{$^{4}$ Universidade Federal de Goi\'as, INF, s/n, 74690-900, Goi\^ania, GO, Brazil.}
\email{priscila.kai@hotmail.com $^{4}$}

\thanks{$^{\ast}$ 
Corresponding author}

\keywords{gradient $k$-Yamabe solitons, Yamabe solitons, $\sigma_{k}$-curvature, $k$-Yamabe problem, invariant solutions, complete examples.}

\subjclass[2010]{53C21, 53C50, 53C25} 

\begin{document}

\begin{abstract}
The purpose of this paper is to study gradient $k$-Yamabe solitons conformal to pseudo-Euclidean space. We characterize all such solitons invariant under the action of an $(n-1)$-dimensional translation group. For rotational invariant solutions, we provide the classification of solitons with null curvatures. As an application, we construct infinitely many explicit examples of  geodesically complete steady gradient $k$-Yamabe solitons conformal to the Lorentzian space.
\end{abstract}
\maketitle
\section{Introduction and main results}
\label{intro}
The concept of gradient $k$-Yamabe soliton, introduced in the celebrated work \cite{catino2012global}, corresponds to a natural generalization of gradient Yamabe solitons. We recall that a pseudo-Riemannian manifold $(M^{n},g)$ is a \textit{gradient $k$-Yamabe soliton} if it admits a constant $\lambda\in\mathbb{R}$ and a function $f\in C^{\infty}(M)$ satisfying the equation 
\begin{equation*}
\nabla^{2}f=2(n-1)(\sigma_{k}-\lambda)g,
\end{equation*}
where $\nabla^{2}f$ and $\sigma_{k}$ stand, respectively, for the Hessian of $f$ and the $\sigma_{k}$-curvature of $g$. Recall that, if we denote by $\lambda_{1},\lambda_{2},\dots,\lambda_{n}$ the
eigenvalues of the symmetric endomorphism $g^{-1}A_{g}$, where $A_{g}$ is the Schouten tensor defined by
\begin{equation*}
A_{g}=\frac{1}{n-2}\left[Ric_{g}-\frac{scal_{g}}{2(n-1)}g\right],
\end{equation*}
then the $\sigma_{k}$-curvature of $g$ is defined as the $k$-th symmetric elementary function of $\lambda_{1},\dots,\lambda_{n}$, namely
\begin{equation*}
\sigma_{k}=\sigma_{k}(g^{-1}A_{g})=\sum_{i_{1}<\dots< i_{k}}\lambda_{i_{1}}\cdot \dots\cdot \lambda_{i_{k}}, \quad \text{for}\quad 1\leq k\leq n.
\end{equation*}
Since $\sigma_{1}$ is the trace of $g^{-1}A_{g}$, the gradient $1$-Yamabe solitons simply correspond to gradient Yamabe solitons \cite{chow1992yamabe,daskalopoulos2013classification, di2008yamabe, hamilton1988ricci, ma2012remarks,tokura2018warped}. As usual, the quadruple $(M^{n}, g, f, \lambda)$ is classified into three types according to the sign of
$\lambda$: expanding if $\lambda<0$, steady if $\lambda= 0$ and shrinking if $\lambda> 0$.  Moreover, when $f$ is a constant function the soliton is called \textit{trivial}.

In recent years, many efforts have been devoted to study the geometry of Riemannian Yamabe solitons and their generalizations. For instance, Hsu in \cite{hsu2012note} shown that any compact gradient $1$-Yamabe soliton is trivial. For $k>1$, the extension of the previous result was recently investigated. Catino \textit{et al.} \cite{catino2012global} proved that any compact
gradient $k$-Yamabe soliton with nonnegative Ricci curvature is trivial. Bo
\textit{et al.} \cite{bo2018k} also proved that any compact gradient $k$-Yamabe
soliton with negative constant scalar curvature necessarily has constant
$\sigma_{k}$-curvature. The previous results were generalized in \cite{tokura2021triviality}, where was shown that any compact gradient $k$-Yamabe
soliton must be trivial.


For the noncompact case, Catino et al. \cite{catino2012global} provide an important relation between gradient $k$-Yamabe solitons and conformally flat spaces; its result establishes that any complete, noncompact gradient $k$–Yamabe soliton with nonnegative Ricci tensor and positive at some point is rotationally symmetric and globally conformally equivalent to Euclidean space.  In the context of conformally flat spaces, Neto and Tenenblat \cite{neto2018gradient} study invariant by translation solutions for gradient $1$-Yamabe solitons. They reduced a system of PDEs, that comes from the corresponding $1$-Yamabe soliton equation, to a system of ODEs by considering a function invariant under translations and, as a result, a complete classification and infinitely many examples are obtained. It is considered this approach in other works like \cite{barbosa2014gradient,barboza2018invariant,batista2019warped, bonfim2019quasi,de2017gradient,leandro2017invariant, leandro2020invariant,tokura2018warped}. In general, the technique of transforming a PDE system into an ODE or a PDE with less independent variables is known as an \textit {ansatz} and an important method for generating ansatz is based in the theory of Lie point symmetry groups for PDE \cite{olver2000applications}.

In this paper, we focus our attention on gradient $k$-Yamabe solitons conformal to pseudo-Euclidean space whose solutions are invariant under the action of an $(n-1)$- dimensional translation group or invariant under the action of a pseudo-orthogonal group. First, we work with the invariant by translation case. More precisely, we consider the pseudo-Riemannian metric
\begin{equation*}
\delta=\sum_{i=1}^{n}\varepsilon_{i}dx_{i}\otimes dx_{i},
\end{equation*}
in coordinates $x=(x_{1},\dots,x_{n})$ of $\mathbb{R}^{n}$, where $n\geq 2$, $\varepsilon_{i}=\pm1$. For an arbitrary choice of non zero vector $\alpha=(\alpha_{1},\dots,\alpha_{n})$ we define the translation function $\xi:\mathbb{R}^{n}\rightarrow\mathbb{R}$ by
\begin{equation*}
\xi(x_{1},\dots,x_{n})=\alpha_{1}x_{1}+\dots+\alpha_{n}x_{n}.
\end{equation*}

Next, we assume that $\mathbb{R}^{n}$ admits a group of symmetries consisting of translations \cite{olver2000applications} and we then look for smooth functions $\varphi, f:(a,b)\subset\mathbb{R}\rightarrow\mathbb{R}$, with $\varphi>0$, such that the compositions 
\[f=f\circ\xi:\xi^{-1}(a,b)\rightarrow\mathbb{R},\qquad \varphi=\varphi\circ\xi: \xi^{-1}(a,b)\rightarrow\mathbb{R},\]
satisfies the gradient $k$-Yamabe soliton equation

\begin{equation}\label{eq fundamental}\nabla^{2}f=2(n-1)\left(\sigma_{k}-\lambda\right)\frac{\delta}{\varphi^2},
\end{equation}
or equivalently, 

\[f_{x_{i},x_{j}}-\sum_{k=1}^{n}\Gamma_{ij}^{k}f_{x_{k}}=2(n-1)\left[\sum_{i_{1}<\dots< i_{k}}\lambda_{i_{1}}\cdot \dots\cdot \lambda_{i_{k}}-\lambda\right]\frac{\delta_{ij}}{\varphi^2},\qquad i,j\in\{1,\dots,n\},\]
where $\varphi_{x_{i}}$, $\varphi_{x_{i},x_{j}}$ denote the derivatives of $\varphi$ with respect the variables $x_i$ and $x_{i}x_{j}$ respectively, and $\Gamma_{ij}^{k}$ correspond to Christoffel symbols on the conformal metric $\delta\varphi^{-2}$. What has been said above is summed up in the next results.

\begin{theorem}\label{translacao}With $(\mathbb{R}^{n},\delta)$ and $f=f\circ\xi$, $ \varphi=\varphi\circ\xi$ as above, the manifold $\xi^{-1}(a,b)\subset \mathbb{R}^n$ furnished with the metric tensor
	\begin{equation*}
	g=\frac{\delta}{\varphi(x)^2},
	\end{equation*}
	is a gradient $k$-Yamabe soliton if, and only if, 
	\begin{equation}\label{eq:01}
	f''+2\frac{f'\varphi'}{\varphi}=0,
	\end{equation}
	\begin{equation}\label{eq:02}
	b_{n,k}\left[k\varphi\varphi''-\frac{n}{2}(\varphi')^{2}\right](\varphi')^{2(k-1)}||\alpha||^{2k}+\frac{\varphi\varphi'f'}{2(n-1)}||\alpha||^{2}=\lambda,
	\end{equation}
	where 
	\begin{equation*}b_{n,k}=\frac{(n-1)!}{k!(n-k)!}(-1)^{k-1}\frac{1}{2^{k-1}}.
	\end{equation*}
\end{theorem}

\begin{observation}The previous theorem can be view as an extension of Theorem 2 of \textup{\cite{neto2018gradient}}, which may be obtained for the particular case $k=1$.
\end{observation}

In the following results we provide the solutions of systems \eqref{eq:01} and \eqref{eq:02} when $\lambda=0$, i.e., in the steady case.

\begin{theorem}\label{theorem 22}With the considerations of Theorem \ref{translacao}, if $\alpha=\sum_{i=1}^{n}\alpha_{i}\frac{\partial}{\partial x_{i}}$ is light-like vector, i.e., $||\alpha||^2=0$, then the gradient $k$-Yamabe soliton is steady. Moreover, given the conformal function $\varphi$, the potential function is given by
	\begin{equation*}
	f(\xi)=\int \frac{c}{\varphi^2(\xi)}d\xi+d,\qquad  c, d\in \mathbb{R}.
	\end{equation*}
\end{theorem}

In the case $||\alpha||^2\neq0$, the behavior of $\varphi$ is classified according to the relation between $2k$
and the dimension $n$. First, we consider the case $n\neq 2k$.
\begin{theorem}\label{theorem 23} Let $(\mathbb{R}^n, \delta_{ij}\varphi^{-2})$, $n\neq2k$, be a steady gradient $k$-Yamabe soliton with $f$ as potential function. Then $\varphi$ and $f$ are invariant under an $(n-1)$-dimensional translation group whose basic invariant is $\xi=\sum\alpha_{i}x_{i}$ and $\alpha=\sum_{i=1}^{n}\alpha_{i}\frac{\partial}{\partial x_{i}}$ is a space-like or time-like vector if, and only if, $\varphi$ and $f$ verify

\begin{itemize}[leftmargin=20pt,align=left,labelwidth=\parindent,labelsep=0pt]
    \item [\textup{(a)}]\hspace{0,1cm} In the case $\varphi'=0$, 
    
	\begin{equation*}
	\varphi(\xi)=b,\qquad f(\xi)=c \xi+d,\qquad   c,d\in \mathbb{R},\qquad b\in(0,\infty).
	\end{equation*}
	
    \item [\textup{(b)}]\hspace{0,1cm}In the case $\varphi'\neq0$,

	\begin{equation}\label{novo3}
	f(\xi)=\int \frac{c}{\varphi^2(\xi)}d\xi+d,\qquad  c,d\in \mathbb{R},
	\end{equation}
	and
\begin{equation}\label{novo4}
	\int\frac{d\varphi}{\varphi^{\frac{n}{2k}}\left[\frac{(n-2k)p(1-2k)}{2nk-n+2k}\varphi^{\frac{2-n}{2}\frac{2nk-n+2k}{n-2k}}+c_{1}\right]^{\frac{1}{2k-1}}}=\frac{2k}{2k-n}\xi+ c_{2},\qquad c_{1}, c_{2}\in\mathbb{R},
	\end{equation}	
where	
\[p=\frac{c}{2(n-1)k b_{n,k}} \left(\frac{2k-n}{2k}\right)^{2(k-1)}\]

\end{itemize}

\end{theorem}

Now, for $n=2k$, we have the following result.

\begin{theorem}\label{theorem 24}Let $(\mathbb{R}^n, \delta_{ij}\varphi^{-2})$, $n=2k$, be a steady gradient $k$-Yamabe soliton with $f$ as potential function. Then $\varphi$ and $f$ are invariant under an $(n-1)$-dimensional translation group whose basic invariant is $\xi=\sum\alpha_{i}x_{i}$ and $\alpha=\sum_{i=1}^{n}\alpha_{i}\frac{\partial}{\partial x_{i}}$ is a space-like or time-like vector if, and only if, $\varphi$ and $f$ satisfy

\begin{itemize}[leftmargin=20pt,align=left,labelwidth=\parindent,labelsep=0pt]
    \item [\textup{(a)}]\hspace{0,1cm} In the case $\varphi'=0$, 
    
	\begin{equation*}
	\varphi(\xi)=c_{1},\qquad f(\xi)=c \xi+d,\qquad   c,d\in \mathbb{R},\qquad c_{1}\in(0,\infty).
	\end{equation*}

    \item [\textup{(b)}]\hspace{0,1cm}In the case $\varphi'\neq0$,

	\begin{equation}\label{21}
	f(\xi)=\int \frac{c}{\varphi^2(\xi)}d\xi+d,\qquad  c,d\in \mathbb{R},
	\end{equation}
	and
		\begin{equation}\label{22}\int\dfrac{d\varphi}{\sqrt[n-1]{\dfrac{c}{b_{n,k}n^2\varphi}+c_{1}\varphi^{n-1}}}=\xi+c_{2},\qquad c_{1},c_{2}\in\mathbb{R}.
		\end{equation}
		
\end{itemize}
\end{theorem}


Following, we deal with the rotational case. In the same way as in the invariant by translation group, we consider $\mathbb{R}^n$ endowed with the pseudo-Euclidean metric $\delta=\sum_{i=1}^{n}\varepsilon_{i}dx_{i}\otimes dx_{i}$ and we define the rotational function $r:\mathbb{R}^{n}\rightarrow\mathbb{R}$ by
\begin{equation*}
 r(x_{1},\dots,x_{n})=\varepsilon_{1}x_{1}^{2}+\dots+\varepsilon_{n}x_{n}^{2}.
\end{equation*}
Next, we look for smooth functions $\varphi, f:(a,b)\subset\mathbb{R}\rightarrow\mathbb{R}$, with $\varphi>0$, such that the compositions 
\[f=f\circ r:r^{-1}(a,b)\rightarrow\mathbb{R},\qquad \varphi=\varphi\circ r: r^{-1}(a,b)\rightarrow\mathbb{R},\]
satisfies the gradient $k$-Yamabe soliton equation \eqref{eq fundamental}.

\begin{theorem}\label{rotacao}Let $(\mathbb{R}^n,\delta_{ij}\varphi^{-2})$ be a conformal to pseudo-Euclidean space. Consider smooth functions $\varphi(r)$, $f(r)$, where $r=\sum \varepsilon_{i}x_{i}^2$. Then $\delta\varphi^{-2}$ is a gradient $k$-Yamabe soliton with $f$ as a potential function if, and only if, $\varphi$ and $f$ satisfy
	\begin{equation}\label{eq:011}
	f''+2\frac{f'\varphi'}{\varphi}=0,
	\end{equation}
	\begin{equation}\label{eq:022}
	c_{n,k}\left[\varphi\varphi'-r(\varphi')^{2}\right]^{k-1}(2n\varphi\varphi'-2nr(\varphi')^{2}+4kr\varphi\varphi'')-\frac{\varphi^{2}f'}{n-1}+\frac{2rf'\varphi'\varphi}{n-1}=\lambda,
	\end{equation}
	where 
	\begin{equation*}c_{n,k}=\frac{(n-1)!}{k!(n-k)!}2^{k-1}.
	\end{equation*}
\end{theorem}

The next result provide the solutions of \eqref{eq:011} and \eqref{eq:022} for gradient $k$-Yamabe solitons with null curvatures.

\begin{theorem}\label{ultimo}Let $(\mathbb{R}^n,\delta_{ij}\varphi^{-2})$ be an invariant by rotation gradient $k$-Yamabe soliton with $f$ as potential function and null $\sigma_{1}$, $\sigma_{s}$, for $s\in\{2,\dots,n\}$. Then $(\mathbb{R}^n, \delta_{ij}\varphi^{-2})$ is  one of the following solitons:

\begin{itemize}
    \item [\textup{(a)}] The Gaussian soliton $(\mathbb{R}^{n},\delta_{ij}\varphi^{-2})$ with potential function and conformal factor given, respectively, by
\[f(r)=\dfrac{(n-1)\lambda}{c_{2}^{2}}r+c_{1}, \qquad \varphi(r)=c_{2},\qquad c_{1}\in\mathbb{R}, \quad c_{2}\in(0,\infty).\]
 
    \item [\textup{(b)}] The soliton $\left(\mathbb{R}^n, \delta_{ij}\varphi^{-2}\right)$ with potential function and conformal factor given, respectively, by
\[f(r)=-\dfrac{(n-1)\lambda}{c_{0}^{2}r}+c_{1},\quad \varphi(r)=c_{0}r,\qquad c_{0}\in(0,\infty), \quad c_{1}\in\mathbb{R}.\]
\end{itemize}
\end{theorem}

We finalize the session with a sufficient condition for the conformal factor $\varphi$ in Theorem \ref{translacao} and Theorem \ref{rotacao} to produce complete metrics in the Riemannian context.

\begin{theorem}\label{final} Let $(\mathbb{R}^n, \delta_{ij}\varphi^{-2})$, be a Riemannian gradient $k$-Yamabe soliton with $f$ as potential function. If $f$ and $\varphi$ are invariant  solutions in Theorem  \ref{translacao} or Theorem \ref{rotacao} with conformal factor satisfying $0<|\varphi|\leq L$ for some constant $L$, then the gradient $k$-Yamabe soliton metric $\delta_{ij}\varphi^{-2} $ is complete.
    
\end{theorem}

\section{Examples}

Before proving our main results, we provide examples illustrating the above theorems. It is worth pointing out that geodesic completeness in pseudo-Riemannian manifolds is an essential concept for studying singularities in general relativity. However, obtaining the geodesics and their singularities explicitly in a pseudo-Riemannian manifold is a difficult task.

Our first example provides geodesically complete steady gradient $k$-Yamabe solitons conformal to the Lorentzian space (see section 3 for more details).


\begin{example}\label{examplecomplete}Consider the Lorentzian space $(\mathbb{R}^{n},g)$ with coordinates $(x_{1},\dots,x_{n})$ and signature $\varepsilon_{1}=-1$, $\varepsilon_{i}=1$ for all $i\in\{2,\dots,n\}$. Let $\xi=x_{1}+x_{2}$ and choose $\theta\in\mathbb{N}$. Then, from Theorem \ref{theorem 22}, the functions 
	\begin{equation*}f(\xi)=c\xi+2c\frac{\xi^{2\theta+1}}{2\theta+1}+c\frac{\xi^{4\theta+1}}{4\theta+1}, \qquad  \varphi(\xi)=\frac{1}{1+\xi^{2\theta}}, \quad c\in\mathbb{R},
	\end{equation*}
defines a family of geodesically complete steady gradient $k$-Yamabe soliton with potential function $f$(see section 3).
\end{example}

\begin{example}In Theorem \ref{theorem 23}, consider $(\mathbb{R}^{n},\delta_{ij}\varphi^{-2})$, $n\neq 2k$, with  $c=0$, then the functions
	\begin{equation*}f(\xi)=c_{0}, \qquad  \varphi(\xi)^{\frac{2k-n}{2k}}=\xi+c_{1},\qquad c_{0}, c_{1}\in\mathbb{R}.
	\end{equation*}
	provide a steady gradient $k$-Yamabe soliton in the semi-space $\xi+c_{1}>0$.
\end{example}

\begin{example}In Theorem \ref{theorem 24}, consider $(\mathbb{R}^{n},\delta_{ij}\varphi^{-2})$ with $n=2k$, $k>1$, $c_{1}=0$ and $c=b_{n,k}n^2$, then the functions
	\begin{equation*}f(\xi)=-\frac{c(n-1)}{n-2}\left(\frac{n}{n-1}\xi+c_{4}\right)^{\frac{2}{n}-1}+c_{3}, \qquad  \varphi(\xi)=\sqrt[\frac{n-1}{n}]{\frac{n}{n-1}\xi+c_{4}},\qquad  c_{3}, c_{4}\in\mathbb{R}.
	\end{equation*}
	provide a steady gradient $k$-Yamabe soliton in the semi-space $n\xi+(n-1)c_{4}>0$.
\end{example}

\begin{example}\label{gauss}\textup{(Gaussian soliton)}The Gaussian soliton on $\mathbb{R}^{n}$ is given by
\[g_{ij}=\delta_{ij},\qquad f(x)=\frac{\lambda}{2}|x|^{2}+c_{1},\qquad c_{1}\in\mathbb{R}.\]
Since $\sigma_{k}(g)=0$ for $k\in\{1,\dots,n\}$, the fundamental equation turns out
\[\nabla^{2}f=\lambda g.\]
It is worth noting that we can get the Gaussian soliton from Theorem \ref{rotacao}. In fact, in  Theorem \ref{rotacao} consider $\varepsilon_{i}=1$ for all $i\in \{1,2,\dots, n\}$, then the functions
	\begin{equation*}f(r)=\frac{\lambda}{2}r, \qquad  \varphi(r)=1,
	\end{equation*}
	provide a complete gradient $k$-Yamabe soliton with soliton constant $\lambda$ and potential function $f$.
    
\end{example}

\begin{example}\label{examplecomplete2}In Theorem \ref{rotacao}, consider $(\mathbb{R}^{n}\setminus\{0\},\delta_{ij}\varphi^{-2})$ with $\varepsilon_{i}=1$, $\forall i\in \{1,2,\dots, n\}$ and $k\neq1$, then the functions
	\begin{equation*}f(r)=-\frac{(n-1)\lambda}{c_{0}^{2}r}, \qquad  \varphi(r)=c_{0}r,\qquad c_{0}\in (0,\infty).
	\end{equation*}
	provide a family of Riemannian gradient $k$-Yamabe soliton with soliton constant $\lambda$.
\end{example}

\begin{example}\label{examplecomplete2}In Theorem \ref{rotacao}, consider $(\mathbb{R}^{n}\setminus\{0\},\delta_{ij}\varphi^{-2})$ with $\varepsilon_{i}=1$, $\forall i\in \{1,2,\dots, n\}$ and $k=1$, then the functions
	\begin{equation*}f(r)=\frac{(n-1)(n+2)}{2}\log(1+r)+c_{0}, \qquad  \varphi(r)=\sqrt{1+r},\qquad c_{0}\in \mathbb{R}.
	\end{equation*}
	provide a family of Riemannian gradient $k$-Yamabe soliton with soliton constant $\lambda=\frac{n-2}{2}$. In the particular case in which $n=2$, this soliton is known as Hamilton’s cigar soliton \cite{hamilton1988ricci}.
\end{example}

\section{Proofs}

\begin{myproof}{Theorem}{\ref{translacao}}

	It is well known that for the conformal metric  $\bar{g}=\varphi^{-2}\delta$, the Ricci curvature is given by \cite{besse2007einstein}:	
\begin{equation}\label{Ricci}Ric_{\bar{g}}=\frac{1}{\varphi^{2}}\Big{\{}(n-2)\varphi \nabla^{2}_{\delta}\varphi+[\varphi\Delta_{\delta}\varphi-(n-1)|\nabla_{\delta}\varphi|^{2}]\delta\Big{\}}.
\end{equation}
So, we easily see that the scalar curvature on conformal metric
	is given by
	\begin{equation}\label{escalar conforme}scal_{\bar{g}}=(n-1)(2\varphi\Delta_{\delta}\varphi-n|\nabla_{\delta}\varphi|^{2}).
	\end{equation}
	
		Now, in order to compute the Schouten Tensor on the conformal geometry $A_{\bar{g}}$ we evoke the expression
	\begin{equation*}
A_{\bar{g}}=\frac{1}{n-2}\left(Ric_{\bar{g}}-\frac{scal_{\bar{g}}}{2(n-1)}\bar{g}\right).
\end{equation*}
Therefore, from \eqref{Ricci} and \eqref{escalar conforme} we deduce that
		\begin{equation*}A_{\bar{g}}=\frac{\nabla^{2}_{\delta}\varphi}{\varphi}-\frac{|\nabla_{\delta}\varphi|^{2}}{2\varphi^{2}}\delta.
	\end{equation*}

Throughout this work we will denote by $\varphi_{x_{i}}$, $\varphi_{x_{i},x_{j}}$ the derivatives of $\varphi$ with respect the variables $x_i$ and $x_{i}x_{j}$, respectively. That being said, since we are assuming that $\varphi(\xi)$ and $f(\xi)$ are functions of $\xi=\alpha_{1}x_{1}+\dots+\alpha_{n}x_{n}$, we get

	\begin{equation*}
	\varphi_{,x_{i}}=\varphi'\alpha_{i},\qquad  f_{,x_{i}}=f'\alpha_{i},\qquad \hspace{0.2cm}\varphi_{,x_{i}x_{j}}=\varphi''\alpha_{i}\alpha_{j},\qquad f_{,x_{i}x_{j}}=f''\alpha_{i}\alpha_{j}.\hspace{0.2cm}
	\end{equation*}
Hence

	\begin{equation*}(\bar{g}^{-1}A_{\bar{g}})_{ij}=\varepsilon_{j }\varphi\varphi''\alpha_{i}\alpha_{j}-\frac{1}{2}(\varphi')^{2}||\alpha||^{2}\delta_{ij}.
	\end{equation*}
The eigenvalues of $\bar{g}^{-1}A_{\bar{g}}$ are $\theta=-\frac{1}{2}(\varphi')^{2}||\alpha||^{2}$ with multiplicity $(n-1)$, and $\mu=(\varphi\varphi''-\frac{1}{2}(\varphi')^{2})||\alpha||^{2}$ with multiplicity $1$. The formula for $\sigma_{k}$ can be found easily by the binomial expansion of $(x-\theta)^{n-1}(x-\mu)$

		\begin{eqnarray}\label{k curvature} \sigma_{k}&=&\frac{(n-1)!}{k!(n-k)!}\left[(n-k)\theta+k\mu\right]\theta^{k-1}\nonumber\\
	&=&\frac{(n-1)!}{k!(n-k)!}(-1)^{k-1}\frac{1}{2^{k-1}}\left[k\varphi\varphi''-\frac{n}{2}(\varphi')^{2}\right](\varphi')^{2(k-1)}||\alpha||^{2k}.
	\end{eqnarray}

	Now, in order to compute the Hessian $\nabla^{2}_{\bar{g}}f$ of $f$ relatively to $\bar{g}$ we evoke the expression
	\begin{equation*}
	(\nabla^{2}_{\bar{g}}f)_{ij}=f_{x_{i},x_{j}}-\sum_{k=1}^{n}\Gamma_{ij}^{k}f_{x_{k}},
	\end{equation*}
	where the Christoffel symbol $\Gamma_{ij}^{k}$ for distinct $i,j,k$ are given by
	\begin{equation}\Gamma_{ij}^{k}=0,\ \Gamma_{ij}^{i}=-\frac{\varphi_{x_{j}}}{\varphi},\ \Gamma_{ii}^{k}=\varepsilon_{i}\varepsilon_{k}\frac{\varphi_{x_{k}}}{\varphi}\;\ \mbox{and}\;\ \Gamma_{ii}^{i}=-\frac{\varphi_{x_{i}}}{\varphi}.\nonumber
	\end{equation}
	Therefore,
		\begin{eqnarray}\label{hessian}(\nabla^{2}_{\bar{g}}f)_{ij}&=&f_{x_{i},x_{j}}+\varphi^{-1}(\varphi_{x_{i}}f_{,x_{j}}+\varphi_{x_{j}}f_{x_{i}})-\delta_{ij}\varepsilon_{i}\sum_{k}\varepsilon_{k}\varphi^{-1}\varphi_{x_{k}}f_{x_{k}}\nonumber\\
	&=&\alpha_{i}\alpha_{j}f''+(2\alpha_{i}\alpha_{j}-\delta_{ij}\varepsilon_{i}||\alpha||^{2})\varphi^{-1}\varphi'f'.
	\end{eqnarray}

Substituting \eqref{k curvature} and \eqref{hessian} into \eqref{eq fundamental} and considering $i\neq j$ we obtain

	\begin{equation*}
	\alpha_{i}\alpha_{j}\left(f''+2\frac{\varphi'f'}{\varphi}\right)=0.
	\end{equation*}
If there exist $i,j$, $i\neq j$ such that $\alpha_{i}\alpha_{j}\neq
	0$, then we get
	\begin{equation*}f''+2\frac{f'\varphi'}{\varphi}=0,
	\end{equation*}
which provides equation \eqref{eq:01}. And for $i=j$, substituting \eqref{k curvature} and \eqref{hessian} into \eqref{eq fundamental} we obtain \eqref{eq:02}. 

Now, we need to consider case $\alpha_{k_{0}}=1$, $\alpha_{k}=0$ for $k\neq k_{0}$. In this case, substituting \eqref{hessian} into \eqref{eq fundamental} we obtain
\begin{equation*}
2(n-1)(\sigma_k-\lambda)\frac{\varepsilon_{i}}{\varphi^2}=-\varepsilon_{i}\frac{\varphi'f'}{\varphi},
\end{equation*}
for $i\neq k_{0}$, that is, $\alpha_{i}=0$, and
\begin{equation*}
2(n-1)(\sigma_k-\lambda)\frac{\varepsilon_{k_0}}{\varphi^2}=f''+(2-\varepsilon_{k_0})\frac{\varphi'f'}{\varphi},
\end{equation*}
for $i=k_{0}$, that is, $\alpha_{k_{0}}=1$. However, this equations are equivalent to \eqref{eq:01} and \eqref{eq:02}. This completes the demonstration.
\end{myproof}

\begin{myproof}{Theorem}{\ref{theorem 22}} Since $||\alpha||^{2}=0$, we have by equation \eqref{eq:02} of Theorem \ref{translacao} that $\lambda=0$. On the other hand, given $\varphi$ we have from \eqref{eq:01} that
\begin{equation*}
    f''+2\frac{\varphi'f'}{\varphi}=0.
\end{equation*}
Integrating we get
\begin{equation*}
f(\xi)=\int\frac{c}{\varphi^{2}(\xi)}d\xi+d,\qquad c,d\in \mathbb{R}.
\end{equation*}
\end{myproof}

\begin{myproof}{Theorem}{\ref{theorem 23}} 
Item (a): Since $\varphi'=0$ we have that \eqref{eq:02} is trivially satisfied and from \eqref{eq:01} we conclude that $f(\xi)=c\xi+d$ for some $c,d\in\mathbb{R}$.

Item (b): From equation \eqref{eq:01} of Theorem \ref{translacao} we deduce that
\begin{equation}\label{novo1}
    f'(\xi)=\frac{c}{\varphi^{2}(\xi)}, \qquad c\in\mathbb{R},
\end{equation}
and then
\begin{equation*}
    f(\xi)=\int\frac{c}{\varphi^{2}(\xi)}d\xi+d,\qquad c,d\in \mathbb{R},
\end{equation*}
which provide equation \eqref{novo3}.

Next, without loss of generality, we may consider $||\alpha||^2=\pm1$. Then, substituting  \eqref{novo1} into \eqref{eq:02} we deduce
\begin{equation}\label{novo5}\left[\varphi\varphi''-\frac{n}{2k}(\varphi')^{2}\right](\varphi')^{2(k-1)}+\frac{c}{2kb_{n,k}(n-1)}\frac{\varphi'}{\varphi}=0.
\end{equation}
Considering $v(\xi)=\varphi(\xi)^{1-\frac{n}{2k}}$, we obtain from \eqref{novo5} the following equivalent condition
\begin{equation}\label{novo6}
    v''+p (v')^{3-2k}v^{\frac{2nk-2n+4k}{n-2k}}=0.
\end{equation}
where 
\[p=\frac{c}{2(n-1)k b_{n,k}} \left(\frac{2k-n}{2k}\right)^{2(k-1)}\]
Now, from one more change $w(v)=(v')^{3-2k}$, we obtain that \eqref{novo6} is equivalent to the following first order differential equation

\begin{equation*}
w'(v)+(3-2k)p w(v)^{\frac{4-4k}{3-2k}}v^{\frac{2nk-2n+4k}{n-2k}}=0,
\end{equation*}
whose solution is
\begin{equation}\label{novo10}
w(v)=\left[\frac{(n-2k)p(1-2k)}{2nk-n+2k}v^{\frac{2nk-n+2k}{n-2k}}+c_{1}\right]^{\frac{2k-3}{1-2k}}, \qquad c_{1}\in\mathbb{R}.
\end{equation}
Replacing \eqref{novo10} back into $w(v)=(v')^{3-2k}$ we deduce that 
\[v'=\left[\frac{(n-2k)p(1-2k)}{2nk-n+2k}v^{\frac{2nk-n+2k}{n-2k}}+c_{1}\right]^{\frac{1}{2k-1}}.\]
This implies that
\[\int\frac{dv}{\left[\frac{(n-2k)p(1-2k)}{2nk-n+2k}v^{\frac{2nk-n+2k}{n-2k}}+c_{1}\right]^{\frac{1}{2k-1}}}=\xi+ c_{2},\qquad c_{2}\in\mathbb{R}.\]
Therefore, it follows from $v(\xi)=\varphi(\xi)^{1-\frac{n}{2k}}$ that

\[\int\frac{d\varphi}{\varphi^{\frac{n}{2k}}\left[\frac{(n-2k)p(1-2k)}{2nk-n+2k}\varphi^{\frac{2-n}{2}\frac{2nk-n+2k}{n-2k}}+c_{1}\right]^{\frac{1}{2k-1}}}=\frac{2k}{2k-n}\xi+ c_{3},\qquad c_{3}\in\mathbb{R}.\]
which provide equation \eqref{novo4}.

\end{myproof}

\begin{myproof}{Theorem}{\ref{theorem 24}} 
Item (a): The proof is analogous to the proof of item (a) in Theorem \ref{theorem 23}.

Item (b):  From equation \eqref{eq:01} of Theorem \ref{translacao} we deduce that
\begin{equation}\label{novo11}
    f'(\xi)=\frac{c}{\varphi^{2}(\xi)}, \qquad c\in\mathbb{R},
\end{equation}
and then
\begin{equation*}
    f(\xi)=\int\frac{c}{\varphi^{2}(\xi)}d\xi+d,\qquad c,d\in \mathbb{R},
\end{equation*}
which provide equation \eqref{21}.

Substituting  \eqref{novo11} into \eqref{eq:02} and considering $n=2k$, we deduce that 
\begin{equation}\label{23}\left[\varphi\varphi''-(\varphi')^{2}\right](\varphi')^{(n-2)}+\frac{\varphi'}{\varphi}\frac{c}{b_{n,k}n(n-1)}=0.
\end{equation}
Considering 
\begin{equation}\label{24}
    w(\varphi)=\varphi',
\end{equation}
we obtain $w'=\frac{\varphi''}{\varphi'}$. So, equation \eqref{23} is equivalent to 
\[(w'(\varphi)w(\varphi)\varphi-w(\varphi)^2)w(\varphi)^{n-2}+\frac{w(\varphi)}{\varphi}\frac{c}{b_{n,k}n(n-1)}=0,\]
whose solution is 
\begin{equation}\label{25}w(\varphi)=\sqrt[n-1]{\frac{c}{b_{n,k}n^2\varphi}+c_{1}\varphi^{n-1}},\qquad c_{1}\in\mathbb{R}.
\end{equation}
Replacing \eqref{25} back into \eqref{24} we deduce that 
\begin{equation*}
   \varphi'=\sqrt[n-1]{\frac{c}{b_{n,k}n^2\varphi}+c_{1}\varphi^{n-1}},\qquad c_{1}\in\mathbb{R}.
\end{equation*}
This implies that
\[\int \frac{d\varphi}{\sqrt[n-1]{\dfrac{c}{b_{n,k}n^2\varphi}+c_{1}\varphi^{n-1}}}=\xi+c_{2},\qquad c_{2}\in\mathbb{R},\]
which provide equation \eqref{22}.
\end{myproof}

\begin{myproof}{Theorem}{\ref{rotacao}}Since we are assuming that $\varphi(r)$ and $f(r)$ are functions of $r$, where $r=\varepsilon_{1}x_{1}^{2}+\dots+\varepsilon_{n}x_{n}^{2}$, we get

\begin{equation}\label{invariante2}
	\begin{aligned}
	\varphi_{x_{i}}&=2\varepsilon_{i}x_{i}\varphi',\qquad & \varphi_{x_{i},x_{j}}&=4\varepsilon_{i}\varepsilon_{j}x_{i}x_{j}\varphi''+2\varepsilon_{i}\delta_{ij}\varphi',\hspace{0.2cm}\\[10pt]
	f_{x_{i}}&=2\varepsilon_{i}x_{i}f',\qquad & f_{x_{i},x_{j}}&=4\varepsilon_{i}\varepsilon_{j}x_{i}x_{j}f''+2\varepsilon_{i}\delta_{ij}f'.\hspace{0.2cm}
	\end{aligned}
	\end{equation}
Substituting \eqref{invariante2} into expression of the Schouten tensor on conformal geometry 

	\begin{equation*}A_{\bar{g}}=\frac{\nabla^{2}_{\delta}\varphi}{\varphi}-\frac{|\nabla_{\delta}\varphi|^{2}}{2\varphi^{2}}\delta,
	\end{equation*}
we deduce that
	\begin{equation*}(\bar{g}^{-1}A_{\bar{g}})_{ij}=4\varepsilon_{j}\varphi\varphi''x_{i}x_{j}+2(\varphi\varphi'-r(\varphi')^{2})\delta_{ij}.
	\end{equation*}
The eigenvalues of $\bar{g}^{-1}A_{\bar{g}}$ are $\lambda=2\varphi\varphi'-2r(\varphi')^{2}$ with multiplicity $(n-1)$, and $\mu=4\varphi\varphi''r+2\varphi\varphi'-2r(\varphi')^{2}$ with multiplicity $1$. The formula for $\sigma_{k}$ can be found easily by the binomial expansion of $(x-\lambda)^{n-1}(x-\mu)$

\begin{equation}\label{k curvature 1}
    \sigma_{k}=\frac{(n-1)!}{k!(n-k)!}2^{k-1}\left[\varphi\varphi'-r(\varphi')^{2}\right]^{k-1}\left[2n\varphi\varphi'-2nr(\varphi')^{2}+4kr\varphi\varphi''\right].
\end{equation}

	Proceeding in a similar way as in the proof of Theorem \ref{translacao}, we obtain the hessian expression on the conformal geometry
		\begin{eqnarray}\label{hessian2}(\nabla^{2}_{\bar{g}}f)_{ij}&=&f_{x_{i},x_{j}}+\varphi^{-1}(\varphi_{x_{i}}f_{,x_{j}}+\varphi_{x_{j}}f_{x_{i}})-\delta_{ij}\varepsilon_{i}\sum_{k}\varepsilon_{k}\varphi^{-1}\varphi_{x_{k}}f_{x_{k}}\nonumber\\
	&=&4\varepsilon_{i}\varepsilon_{j}x_{i}x_{j}\left(f''+2\frac{f'\varphi'}{\varphi}\right)+2\varepsilon_{i}\delta_{ij}\left(f'-2\frac{\varphi'f'}{\varphi}r\right).
	\end{eqnarray}
Substituting \eqref{k curvature 1} and \eqref{hessian2} into \eqref{eq fundamental} and considering $i\neq j$ we obtain
	\begin{equation*}f''+2\frac{f'\varphi'}{\varphi}=0,
	\end{equation*}
which provides equation \eqref{eq:011}, and for $i=j$, substituting \eqref{k curvature 1} and \eqref{hessian2} into \eqref{eq fundamental} we obtain \eqref{eq:022}.

\end{myproof}

\begin{myproof}{Theorem}{\ref{ultimo}} 
The hypothesis $\sigma_{1}=\sigma_{s}=0$, $s\in\{2,\dots,n\}$ implies that

\begin{equation}\label{31} \sigma_{1}=\left[2n\varphi\varphi'-2nr(\varphi')^{2}+4r\varphi\varphi''\right]=0,
\end{equation}
and 
\begin{equation}\label{32}\sigma_{s}=\frac{(n-1)!}{s!(n-s)!}2^{s-1}\left[\varphi\varphi'-r(\varphi')^{2}\right]^{s-1}\left[2n\varphi\varphi'-2nr(\varphi')^{2}+4sr\varphi\varphi''\right]=0.
\end{equation}
	
Therefore, from \eqref{31} and \eqref{32}, we conclude that

\[\left[\varphi\varphi'-r(\varphi')^{2}\right]^{s-1}\left[4(s-1)r\varphi\varphi''\right]=0.\]
whose general solution is given by $\varphi(r)=ar+b$, $a,b\in\mathbb{R}$. We claim that $a=0$ or $b=0$. In fact, replacing $\varphi(r)=ar+b$ into \eqref{31}, we deduce that 
\[2nab=2n(ar+b)a-2nra^2=0\]
which proves the assertion.

Now, suppose that $a=0$, then the conformal factor is constant $\varphi(r)=b>0$. Hence, from \eqref{eq:011} and \eqref{eq:022} we conclude that \[f(r)=\dfrac{(n-1)\lambda}{b^{2}}r+c_{1}, \qquad c_{1}\in\mathbb{R},\]
which provide the proof of item (a).

On the other hand, suppose that $a\neq0$, then $b=0$ and the conformal factor is given by $\varphi(r)=ar$, $a\in(0,\infty)$. Hence, from \eqref{eq:011} and \eqref{eq:022} we conclude that 
\[f(r)=-\dfrac{(n-1)\lambda}{a^{2}r}+c_{1},\qquad c_{1}\in\mathbb{R},\]
which provide the proof of item (b).

\end{myproof}

\begin{myproof}{Theorem}{\ref{final}} Let $(\mathbb{R}^{n}, \bar{g}=\delta\varphi^{-2})$ be the gradient $k$-Yamabe soliton with invariant solutions $f$ and $\varphi$. Since $0<|\varphi|\leq L$, we have that 
\[0<N\leq\frac{1}{L^2}\leq\frac{1}{\varphi^2},\]
for some $N\in(0,\infty)$. This implies that $|v|_{\bar{g}}\geq |v|_{\delta}$ for any vector $v\in \mathbb{R}^n$. Since $(\mathbb{R}^{n},\delta)$ is complete, it follows that $(\mathbb{R}^{n},\delta\varphi^{-2})$ is complete.

\end{myproof}

\begin{myproof1}{}{\ref{examplecomplete}}
Let $(\mathbb{R}^{n},\delta)$ be the standard pseudo-Euclidean space where $\delta=-dx_{1}^{2}+\sum_{i=2}^{n}dx_{i}^{2}$. Take $\theta\in\mathbb{N}$ and consider the functions 	
\begin{equation*} f(\xi)=c\xi+2c\frac{\xi^{2\theta+1}}{2\theta+1}+c\frac{\xi^{4\theta+1}}{4\theta+1}, \qquad  \varphi(\xi)=\frac{1}{1+\xi^{2\theta}},\quad \xi=x_{1}+x_{2},\quad \theta\in\mathbb{N}, \quad c\in\mathbb{R}.
	\end{equation*}
We will now prove that $(\mathbb{R}^{n},\delta\varphi^{-2})$ is geodesically complete by showing that any geodesic
$\gamma(t)=(x_{1}(t),x_{2}(t),\dots, x_{n}(t))\in \mathbb{R}^{n}$ is defined for all $t\in \mathbb{R}$. From the fundamental geodesic equations 
\begin{equation*}
x_{l}''(t)=-\sum_{i,j}\Gamma_{ij}^{l}x_{i}'(t)x_{j}'(t)
\end{equation*}
and the Christoffel symbol $\Gamma_{ij}^{l}$ on the conformal metric
	\begin{equation}\Gamma_{ij}^{l}=0,\ \Gamma_{ij}^{i}=-\frac{\varphi_{x_{j}}}{\varphi},\ \Gamma_{ii}^{l}=\varepsilon_{i}\varepsilon_{l}\frac{\varphi_{x_{l}}}{\varphi}\;\ \mbox{and}\;\ \Gamma_{ii}^{i}=-\frac{\varphi_{x_{i}}}{\varphi},\nonumber
	\end{equation}
we deduce, for $l=1,2,\dots,n$, that

\begin{equation}\label{comp1}
x_{l}''(t)=\frac{1}{\varphi}\left[2(\varphi\circ\gamma)'(t)x_{l}'(t)-\varepsilon_{l}\alpha_{l}\varphi'\circ\xi\circ\gamma(t)\sum_{i=1}^{n}\varepsilon_{i}\left(x_{i}'(t)\right)^{2}\right].
\end{equation}
Therefore, since $\varepsilon_{1}=-1$, $\varepsilon_{i}=1$, $i\geq2$, $\alpha_{1}=\alpha_{2}=1$ and $\alpha_{i}=0$, $i\geq3$, we get

\begin{equation*}
x_{l}''(t)=\frac{1}{\varphi}\left[2(\varphi\circ\gamma)'(t)x_{l}'(t)\right]=-\frac{4\theta(\xi\circ\gamma(t))^{2\theta-1}}{1+(\xi\circ\gamma(t))^{2\theta}}(\xi\circ\gamma)'(t)x_{l}'(t),\quad l\geq3.
\end{equation*}
This implies that
\begin{equation*}
    \left[x_{l}'(t)(1+(\xi\circ\gamma(t))^{2\theta})^{2}\right]'=0,
\end{equation*}
and hence

\begin{equation}\label{sis1}
    x_{l}'(t)=\frac{c_{l}}{(1+(\xi\circ\gamma(t))^{2\theta})^{2}},\qquad c_{l}\in\mathbb{R},\quad l\geq3.
\end{equation}

Now, for $l=1$, we have from \eqref{comp1} that
\begin{eqnarray}\label{sis2}
x_{1}''(t)&=&\frac{1}{\varphi}\left[2(\varphi\circ\gamma)'(t)x_{1}'(t)-\varepsilon_{1}\alpha_{1}\varphi'\circ\xi\circ\gamma(t)\sum_{i=1}^{n}\varepsilon_{i}\left(x_{i}'(t)\right)^{2}\right].\nonumber\\
&=&-\frac{2\theta(\xi\circ\gamma(t))^{2\theta-1}}{1+(\xi\circ\gamma(t))^{2\theta}}\left(2x_{1}'(t)(\xi\circ\gamma)'(t)-(x_{1}'(t))^{2}+(x_{2}'(t))^{2}\right)+\nonumber\\
&&-\frac{2\theta(\xi\circ\gamma(t))^{2\theta-1}}{1+(\xi\circ\gamma(t))^{2\theta}}\frac{1}{(1+(\xi\circ\gamma(t))^{2\theta})^{4}}\sum_{i=3}^{n}c_{i}\nonumber\\
&=&-\frac{2\theta(\xi\circ\gamma(t))^{2\theta-1}}{1+(\xi\circ\gamma(t))^{2\theta}}[(\xi\circ\gamma)'(t)]^{2}-\frac{2\theta(\xi\circ\gamma(t))^{2\theta-1}}{(1+(\xi\circ\gamma(t))^{2\theta})^{5}}\sum_{i=3}^{n}c_{i}.
\end{eqnarray}
Similarly, for $l=2$, we get
\begin{equation}\label{sis3}
    x_{2}''(t)=-\frac{2\theta(\xi\circ\gamma(t))^{2\theta-1}}{1+(\xi\circ\gamma(t))^{2\theta}}[(\xi\circ\gamma)'(t)]^{2}+\frac{2\theta(\xi\circ\gamma(t))^{2\theta-1}}{(1+(\xi\circ\gamma(t))^{2\theta})^{5}}\sum_{i=3}^{n}c_{i}.
\end{equation}
Therefore, 
\begin{equation*}
    (\xi\circ\gamma)''(t)=x_{1}''(t)+x_{2}''(t)=-\frac{4\theta(\xi\circ\gamma(t))^{2\theta-1}}{1+(\xi\circ\gamma(t))^{2\theta}}[(\xi\circ\gamma)'(t)]^{2},
\end{equation*}
which implies that
\begin{equation}\label{sis4}
(\xi\circ\gamma)'(t)=\frac{k_{1}}{(1+(\xi\circ\gamma(t))^{2\theta})^{2}}.
\end{equation}

It follows from \eqref{sis1}, \eqref{sis2}, \eqref{sis3} and \eqref{sis4} that we are looking for the solutions of the system
\begin{equation*}
\left\{\begin{array}{r@{\mskip\thickmuskip}l}
	x_{1}'(t)&= y_{1}(t),\\ [10pt]
	x_{2}'(t)&= y_{2}(t),\\[10pt]
	y_{1}'(t)&=-\dfrac{2\theta (x_{1}(t)+x_{2}(t))^{2\theta-1}}{[1+(x_{1}(t)+x_{2}(t))^{2\theta}]^{5}}(k_{1}^{2}+\sum_{i=3}^{n}c_{i}),\\[10pt]
	y_{2}'(t)&=\dfrac{2\theta}{[1+(x_{1}(t)+x_{2}(t))^{2\theta}]^{2}}(-k_{1}^{2}+\sum_{i=3}^{n}c_{i}),\\[10pt]
	x_{l}'(t)&=\dfrac{c_{l}}{(1+(\xi\circ\gamma(t))^{2\theta})^{2}},\hspace{3cm} l\geq3
	
\end{array} \right.
\end{equation*}

Since the functions
\begin{equation*}
    p(x_{1},x_{2})=\frac{(x_{1}+x_{2})^{2\theta-1}}{[1+(x_{1}+x_{2})^{2\theta}]^{5}},\qquad q(x_{1},x_{2})=\frac{1}{[1+(x_{1}+x_{2})^{2\theta}]^{2}},
\end{equation*}
are of bounded derivative, we conclude that $p$ and $q$ are Lipschitz. Therefore, the solutions of above system exists for all $t\in \mathbb{R}$ and hence all geodesic $\gamma(t)$ are defined for the entire real line, which
means that $(\mathbb{R}^{n},\varphi^{-2}\delta)$ is geodesically complete.

\end{myproof1}

\end{document}